\documentclass{amsart}
\newtheorem{thm}{Theorem}[section]

\newtheorem{lem}[thm]{Lemma}

\theoremstyle{definition}

\theoremstyle{remark}
\newtheorem{rem}[thm]{Remark}
\numberwithin{equation}{section}
\usepackage{graphicx}
\usepackage{amsmath,amsfonts,amssymb}
\usepackage{enumerate}

\newcommand{\bolo}{\textordmasculine}

\begin{document}

\title[The Generalization of Miquel's Theorem]{The Generalization of Miquel's Theorem}%
\author{Anderson R. Vargas}%
\email{}%

\thanks{Author's e-mail address: avargas@mat.puc-rio.br\\
Anderson R. Vargas is a teacher in Colégio Pedro II and a Ph.D. student in Pontifical Catholic University of Rio de Janeiro, with the support of CNPq.}%

\begin{abstract}
This papper aims to present and demonstrate Clifford's version for a generalization of Miquel's theorem with the use of Euclidean geometry arguments only.
\end{abstract}
\maketitle

\section{\textbf{Introduction}}

At the end of his article, Clifford \cite{Clifford1882a} gives some developments that generalize the three circles version of Miquel's theorem and he does give a synthetic proof to this generalization using arguments of projective geometry. The series of propositions given by Clifford are in the following theorem:

\begin{thm}\label{Teo-generalizacao}${}$
\begin{enumerate}[(i)]
\item Given three straight lines, a circle may be drawn through their intersections.
\item Given four straight lines, the four circles so determined meet in a point.
\item Given five straight lines, the five points so found lie on a circle.
\item Given six straight lines, the six circles so determined meet in a point.
\end{enumerate}

That can keep going on indefinitely, that is, if $n\geq2$, $2n$ straight lines determine $2n$ circles all meeting in a point, and for $2n+1$ straight lines the $2n+1$ points so found lie on the same circle.
\end{thm}

\begin{rem}\label{line-remark}
Note that in the set of given straight lines, there is neither a pair of parallel straight lines nor a subset with three straight lines that intersect in one point. That is being considered all along the work, without further ado.
\end{rem}

In order to prove this generalization, we are going to use some theorems proposed by Miquel \cite{Miquel1838-2} and some basic lemmas about a bunch of circles and their intersections, and we will follow the idea proposed by Lebesgue\cite{Lebesgue} in a proof by induction.

\section{\textbf{Preliminaries}}

\begin{figure}[!h]
\center
  \includegraphics*[width=0.6\linewidth]{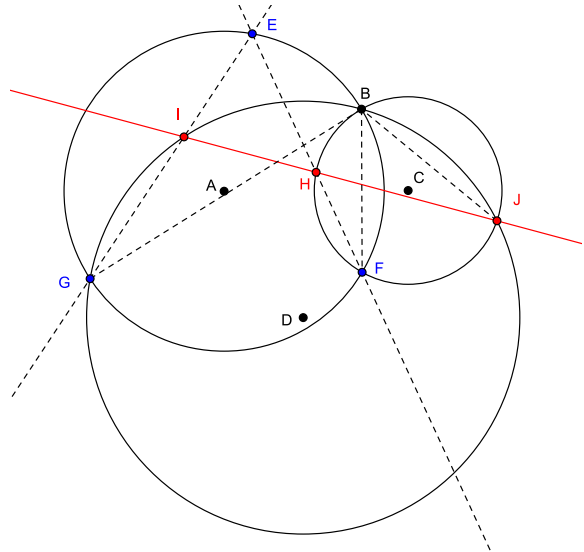}
  \caption{\label{Fig.Miquel1}Miquel's First Theorem}
\end{figure}

\begin{thm}[\textbf{Miquel's First Theorem \cite{Miquel1838-2}}]\label{Miquel1}
Let us consider the circles $A$, $C$ and $D$ (denoted by their centers) which meet in the point $B$. Let us take the point $E$ lying on $A$ and let  $F$ and $G$ be the intersection point between $A$ and $C$, $A$ and $D$, respectively, distinct from $B$. Let us consider the points $H$ and $I$ as the intersections between the straight lines $EF$ and $EG$ and the circles $C$ and $D$, respectively. If $J$ is the intersection point between the circles $C$ and $D$, then the points $H$, $I$ and $J$ are collinear (see Fig.\ref{Fig.Miquel1}).
\end{thm}

\begin{thm}[\textbf{First Reciprocal}]\label{Miquel2}
Let us consider the circles $A$, $C$ and $D$ which meet in the point $B$. Let $J$ be the intersection point between the circles $C$ and $D$. Let us take the points $H$ and $I$ lying on the circles $C$ and $D$, respectively, so that $H$, $I$ and $J$ are collinear. If the points $F$ and $G$ are the intersection points between $A$ and $C$, $A$ and $D$, respectively, distinct from $B$, then the straight lines $FH$ and $GI$ meet in the point $E$ which lies on the circle $A$ (see Fig.\ref{Fig.Miquel1}).
\end{thm}

\begin{figure}[!h]
\center
  \includegraphics*[width=0.6\linewidth]{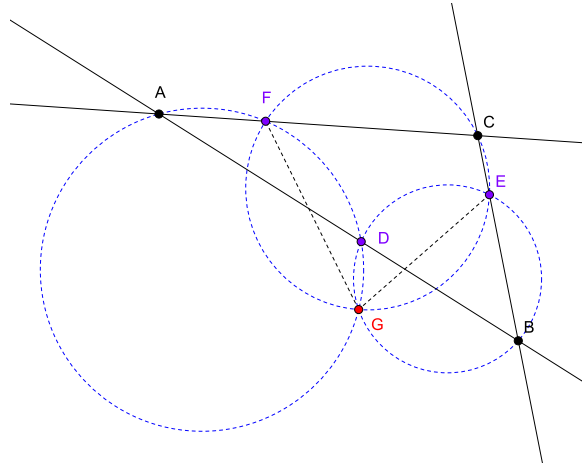}
  \caption{\label{Fig.Miquel3} Miquel's Theorem - second reciprocal}
\end{figure}

\begin{thm}[\textbf{Second Reciprocal}]\label{Miquel3}
Let us consider the points $D$, $E$ and $F$, lying on the straight lines $AB$, $BC$ and $CA$, respectively, so that all points are distinct. Then the three circles $ADF$, $BDE$ and $CEF$ meet in a point $G$ (see Fig.\ref{Fig.Miquel3}).
\end{thm}

The proof of these three Theorems is elementary and can be found in \cite{Miquel1838-2} or \cite{Vargas}.

\begin{figure}[!h]
\center
  \includegraphics*[width=0.6\linewidth]{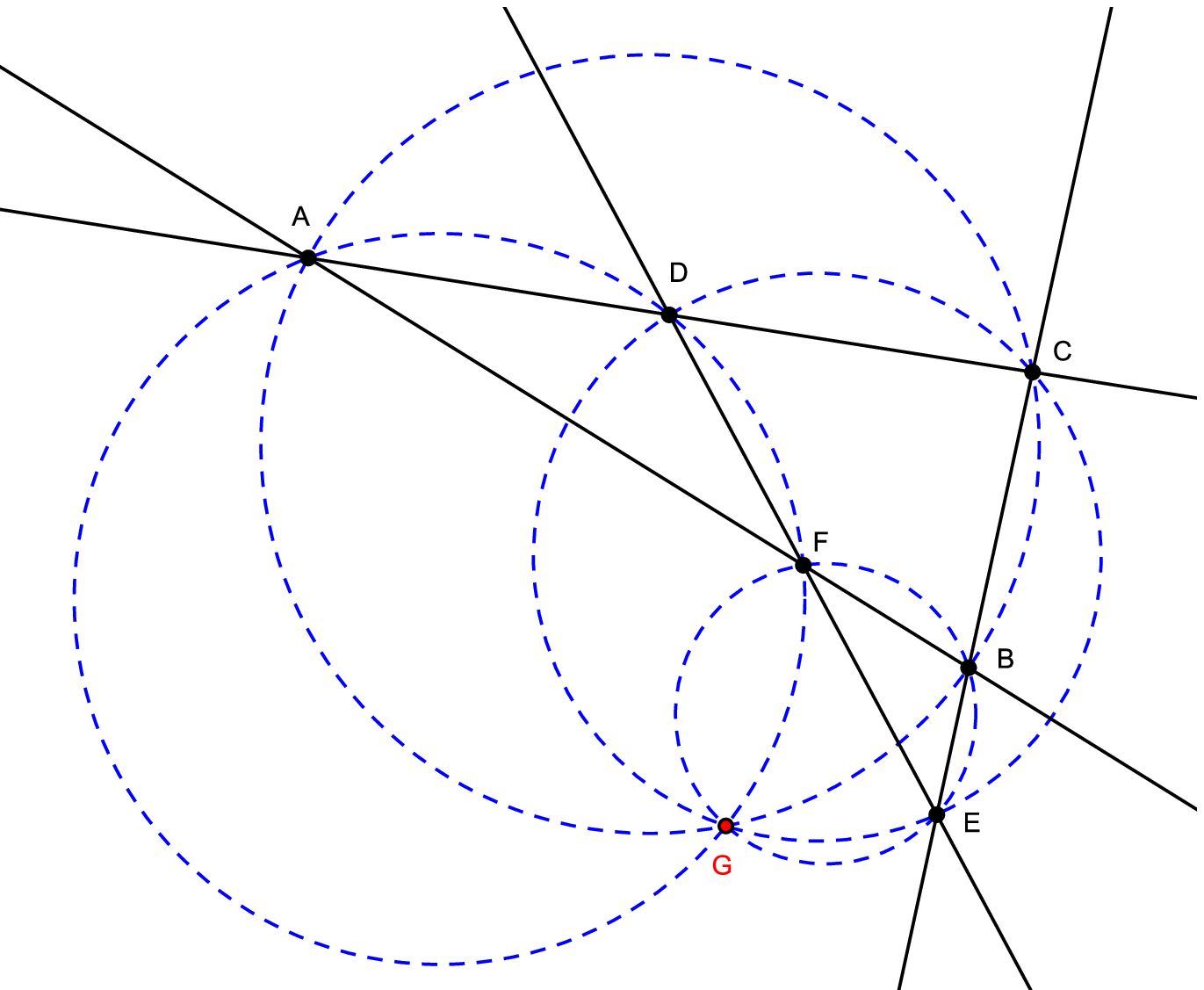}
  \caption{\label{Fig.Miquel4}Theorem \ref{Miquel4}}
\end{figure}

\begin{thm} \label{Miquel4}
Given four straight lines that respect Remark \ref{line-remark}, we have formed four triangles whose vertices are the intersection points between the given lines. Each triangle is inscribable to a circle, then the four circles so determined meet in a point (see Fig.\ref{Fig.Miquel4}).
\end{thm}

$\hspace{-0.55cm}$ \emph{Proof.}
It follows directly by Theorem \ref{Miquel3} applied to triangles $ABC$ and $CDE$.\hfill $\blacksquare$\\

\begin{figure}[!h]
\center
\includegraphics*[width=1\linewidth]{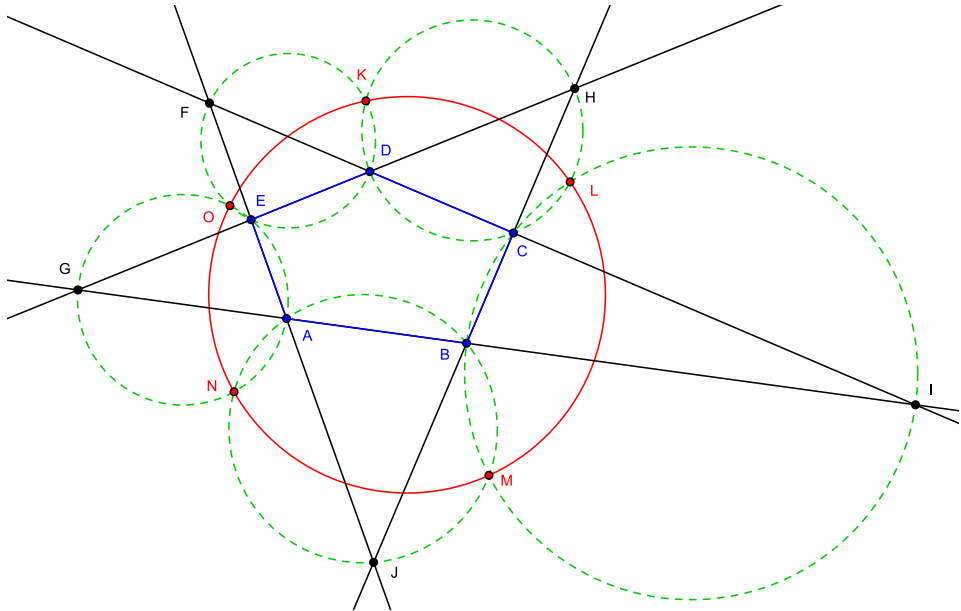}
\caption{\label{Fig.Miquel5}Miquel's Theorem for the pentagon}
\end{figure}

$\vspace{-1cm}$ ${}$\\
\begin{thm}[Miquel's Theorem for the Pentagon]\label{Miquel5}${}$\\
Let $ABCDE$ be a pentagon and let $F$, $G$, $H$, $I$ and $J$, the intersection points between the lines on which lie the pentagon edges. In such a manner we have formed the triangles $ABJ$, $BCI$, $CDH$, $DEF$ and $AEG$, and consequently, the circles circumscribable to them. Then the points $K$, $L$, $M$, $N$ and $O$, which are the intersections points between two adjacent circles, other than pentagon vertices, lie on a circle (see Fig.\ref{Fig.Miquel5}).
\end{thm}

\newpage

\begin{figure}[!h]
    \center
    \includegraphics*[width=\linewidth]{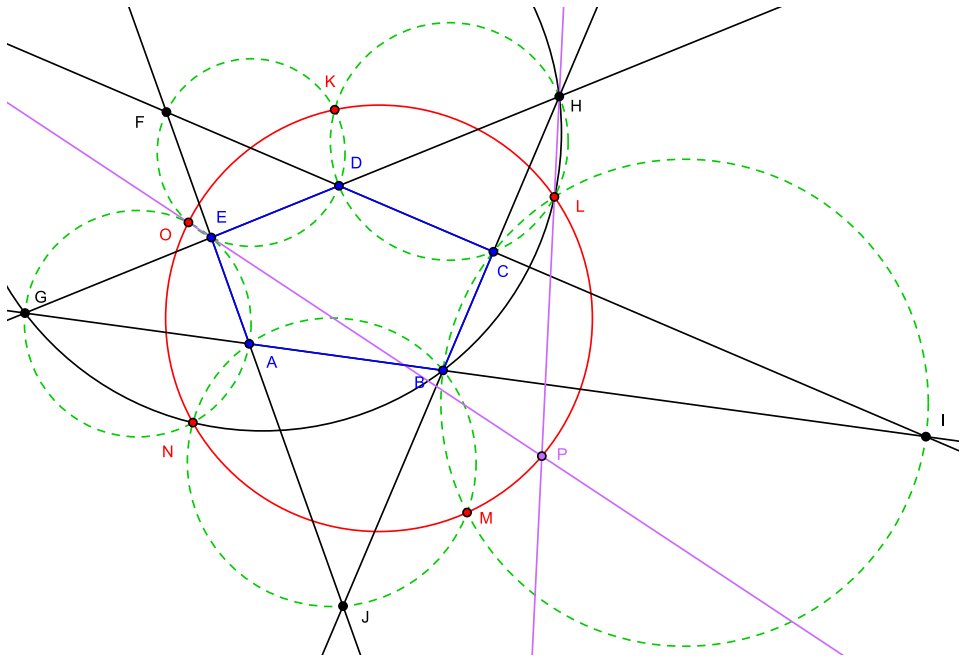}
    \caption{\label{Fig.Miquel6}Miquel's Theorem for the pentagon - proof}
    \end{figure}

\newpage
$\hspace{-0.55cm}$ \emph{Proof.}
Let us consider the circle determined by the points $N$, $L$, $O$. We would like to prove that $K$ lies on it.

Let us take a look on the circle $GBH$ and the quadrilaterals  $GHJA$ and $HGIC$. Then Theorem \ref{Miquel4} implies that this circle contains the points $N$ and $L$ (see Fig.\ref{Fig.Miquel6}).

The circles $NLH$, $NLO$ and $NEO$, meet in the point $N$. The straight line $EH$ contains the point $G$, intersection point between $NLH$ and $NEO$ distinct from $N$. Then Theorem \ref{Miquel2} implies that the lines $OE$ and $HL$ intersect each other in the point $P$, which lies on the circle $NLO$ (see Fig.\ref{Fig.Miquel6}).

Note that the points $O$, $L$, $D$, lie on the straight lines $EP$, $PH$, $HE$, respectively. Then Theorem \ref{Miquel3} implies that the circles $POL$, $EDO$, $HLD$, meet in a point. Consequently, the circle $NLO$ contains the point $K$, which is the intersection point between the circles $EDF$ and $CDH$.

Analogously, we can prove that $NLO$ contains the point $M$.

So we have proved that $K$, $L$, $M$, $N$, $O$, lie on a circle.\hfill $\blacksquare$\\

\begin{figure}[!h]
\center
  \includegraphics*[width=0.9\linewidth]{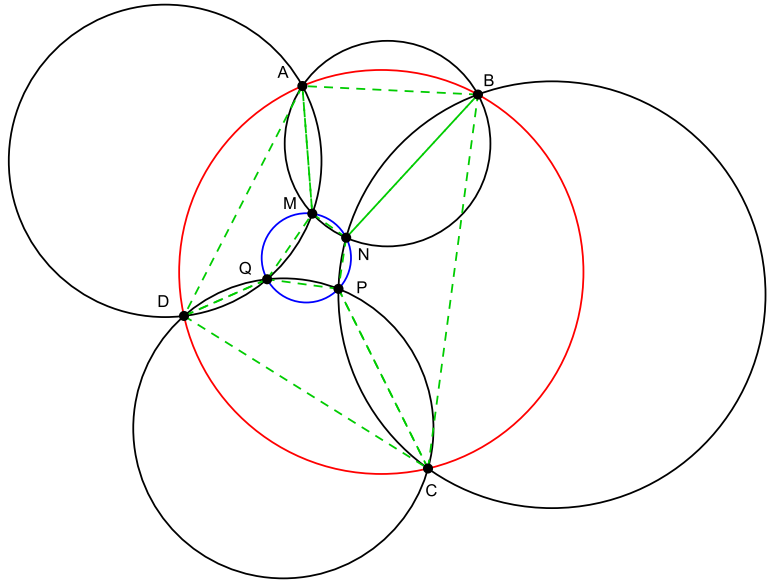}
  \caption{\label{Fig.Lema-4circulos}Four circles Lemma}
\end{figure}

\begin{lem}\label{Lebesgue}
Given four circles $C_1$, $C_2$, $C_3$, $C_4$, let $A$, $B$, $C$, $D$ and $M$, $N$, $P$, $Q$ be the intersection points between $C_i$ and $C_j$, where $i\neq j$ and $i,j\in\{1,2,3,4\}$. If the points $A$, $B$, $C$, $D$ lie on a circle then the points $M$, $N$, $P$, $Q$ also lie on a circle (see Fig.\ref{Fig.Lema-4circulos}).
\end{lem}

\newpage
$\hspace{-0.55cm}$ \emph{Proof.}
Let $A$, $B$, $C$, $D$ be the intersection points lying on a circle (see Fig.\ref{Fig.Lema-4circulos}). That means the quadrilateral $ABCD$ is inscribable, \emph{i.e.},

\begin{equation}\label{eq7}
     m(D\hat{A}B)+m(B\hat{C}D)=180\bolo
    \end{equation}

Let $M$, $N$, $P$, $Q$, be the other intersection points, as it can be seen on Fig.\ref{Fig.Lema-4circulos}. Therefore, the quadrilaterals $AMNB$, $BNPC$, $CPQD$, $DQMA$ are all inscribable and it follows directly that:
    \vspace{-0.05cm}
    \begin{equation}\label{eq8}
     m(M\hat{A}B)+m(B\hat{N}M)=180\bolo
    \end{equation}
    \vspace{-0.7cm}
    \begin{equation}\label{eq9}
     m(D\hat{A}M)+m(M\hat{Q}D)=180\bolo
    \end{equation}
    \vspace{-0.7cm}
    \begin{equation}\label{eq10}
     m(B\hat{C}P)+m(P\hat{N}B)=180\bolo
    \end{equation}
    \vspace{-0.7cm}
    \begin{equation}\label{eq11}
     m(P\hat{C}D)+m(D\hat{Q}P)=180\bolo
    \end{equation}
    \vspace{-0.2cm}

    Let us add the respective members of equations \ref{eq8} to \ref{eq11}, then we have:
    \begin{equation}\label{eq12}
    m(D\hat{A}B)+m(B\hat{C}D)+m(P\hat{N}M)+m(M\hat{Q}P)=720\bolo
    \end{equation}

    The equations \ref{eq7} and \ref{eq12} imply that:
    $$m(P\hat{N}M)+m(M\hat{Q}P)=540\bolo$$

     Observe that we are using a notation for angles that preserve the counterclockwise orientation, which means $M\hat{N}P$ is different from $P\hat{N}M$, but $m(M\hat{N}P)+m(P\hat{N}M)=360\bolo$. Then,
    $$m(M\hat{N}P)+m(P\hat{Q}M)=360\bolo - m(P\hat{N}M)+ 360\bolo - m(M\hat{Q}P)$$
    $$\phantom{m(M\hat{N}P)+m(P\hat{Q}M)}=720\bolo - (m(P\hat{N}M)+ m(M\hat{Q}P))\phantom{aa\;aa}$$
    $$\phantom{m(M\hat{N}P)+m(P\hat{Q}M)}=720\bolo - 540\bolo = 180\bolo\phantom{aaaaaaaaaaaaaaa}$$

    Therefore, $MNPQ$ is inscribable, that is, $M$, $N$, $P$, $Q$, lie on a circle. \hfill $\blacksquare$\\

    The next Lemma can be seen as a particular case of the previous one by considering a line as a circle with infinity radius. But our intention here is to present a readable work for those who are not familiarized with this concept. So, instead of circular points, let us take by hypothesis collinear intersection points, as follows:

\begin{lem}\label{Lebesgue2}
Given four circles $C_1$, $C_2$, $C_3$, $C_4$, let $A$, $B$, $C$, $D$ and $M$, $N$, $P$, $Q$ be the intersection points between $C_i$ and $C_j$, where $i\neq j$ and $i,j\in\{1,2,3,4\}$. If the points $A$, $B$, $C$, $D$ are collinear then the points $M$, $N$, $P$, $Q$ lie on a circle.
\end{lem}

\newpage
\begin{figure}[!h]
\center
  \includegraphics*[width=0.7\linewidth]{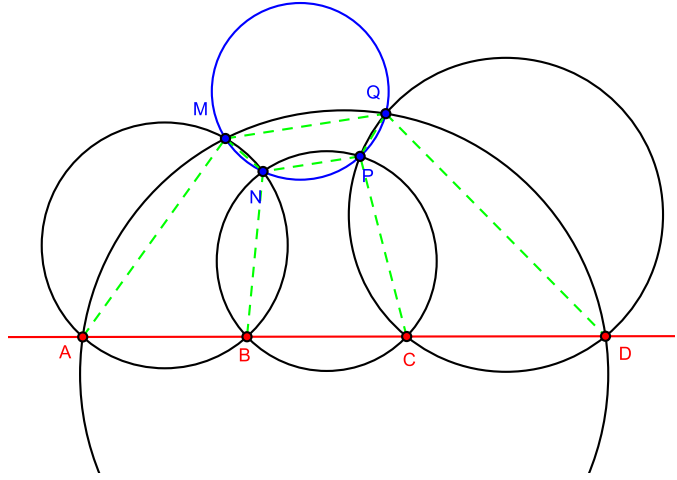}
  \caption{\label{Fig.Lema2}Four circles Lemma with collinearity}
\end{figure}

$\hspace{-0.55cm}$ \emph{Proof.}
     Let $A$, $B$, $C$, $D$, be the collinear intersection points and let $M$, $N$, $P$, $Q$, be the other intersections (see Fig.\ref{Fig.Lema2}). Note that the quadrilaterals $ABNM$, $BCPN$, $CDQP$, $ADQM$, are all inscribable. So the following equations are equivalent,

    \vspace{0.1cm}
    $m(A\hat{M}Q)+m(Q\hat{D}A)=180\bolo$

    \vspace{0.1cm}
    $m(A\hat{M}N)+m(N\hat{M}Q)+m(Q\hat{D}C)=180\bolo$

    \vspace{0.1cm}
    $m(A\hat{M}N)+m(N\hat{M}Q)+180\bolo\,-m(C\hat{P}Q)=180\bolo$

    \vspace{0.1cm}
    $m(A\hat{M}N)+m(N\hat{M}Q)-(360\bolo\,-m(Q\hat{P}C))=0$

    \vspace{0.1cm}
    $m(A\hat{M}N)+m(N\hat{M}Q)+m(Q\hat{P}C)=360\bolo$

    \vspace{0.1cm}
    $m(A\hat{M}N)+m(N\hat{M}Q)+m(Q\hat{P}N)+m(N\hat{P}C)=360\bolo$

    \vspace{0.1cm}
    $m(A\hat{M}N)+m(N\hat{M}Q)+m(Q\hat{P}N)+180\bolo\,-m(C\hat{B}N)=360\bolo$

    \vspace{0.1cm}
    $m(A\hat{M}N)+m(N\hat{M}Q)+m(Q\hat{P}N)-m(C\hat{B}N)=180\bolo$

    \vspace{0.1cm}
    $m(A\hat{M}N)+m(N\hat{M}Q)+m(Q\hat{P}N)-(180\bolo\,-m(N\hat{B}A)=180\bolo$

    \vspace{0.1cm}
    $m(N\hat{M}Q)+m(Q\hat{P}N)+m(A\hat{M}N)+m(N\hat{B}A)-180\bolo=180\bolo$

    \vspace{0.1cm}
    $m(N\hat{M}Q)+m(Q\hat{P}N)=180\bolo$

    \vspace{0.1cm}
    Therefore, $MNPQ$ is inscribable, that is, $M$, $N$, $P$, $Q$, lie on a circle. \hfill $\blacksquare$\\

\section{\textbf{A proof by induction}}

We are going to prove the Theorem \ref{Teo-generalizacao} using an inductive argument over the number of straight lines. Note that the theorem has a thesis for an even number of straight lines and another one for an odd number of straight lines, \emph{i.e.}, if we have $2n$ or $2n+1$ straight lines. In order to fix the notation we are going to do the demonstration for the four items described in the introduction, in this way we will have already proved the first step of the induction.\\

   (i) The first case is well known in Euclidean geometry, three points lie on a circle and there is no need to prove it. However, let us begin fixating our notation.

\begin{figure}[!h]
\center
  \includegraphics*[width=0.6\linewidth]{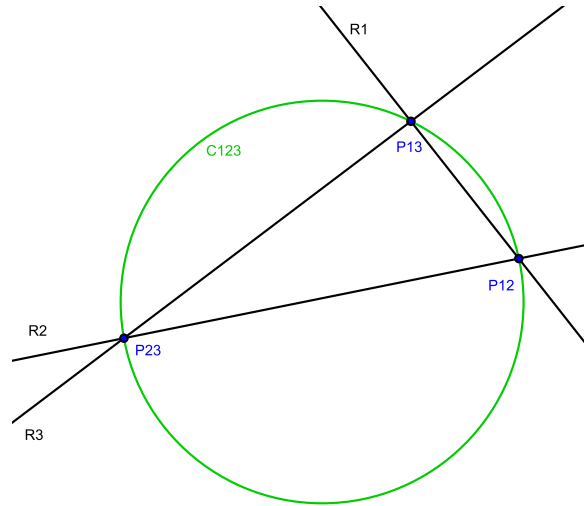}
  \caption{\label{Generalizacao3retas}Three straight lines}
\end{figure}

    Given three straight lines $R_1$, $R_2$ and $R_3$, let $P_{12}$, $P_{13}$ and $P_{23}$ be the intersection points accordingly to the indexes (see Fig.\ref{Generalizacao3retas}). So these points lie on the circle $C_{123}$.\\

\begin{figure}[!h]
\center
  \includegraphics*[width=0.8\linewidth]{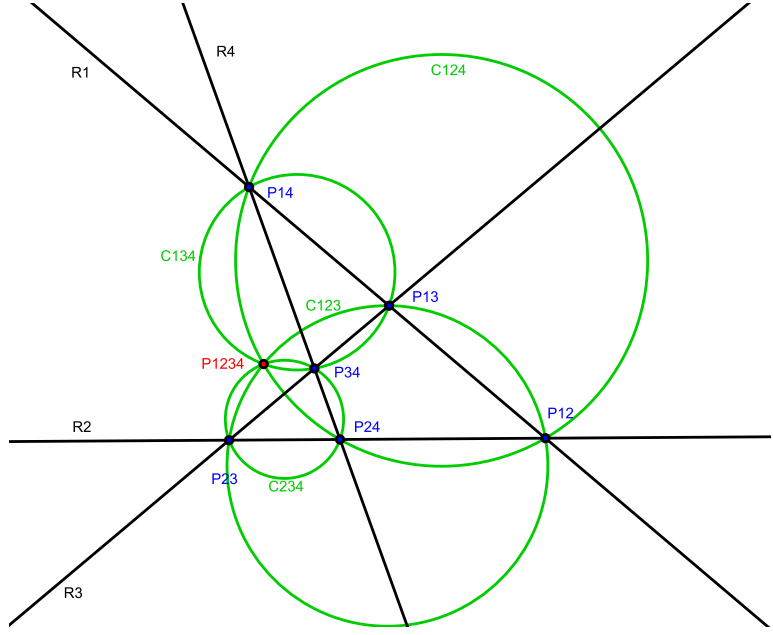}
  \caption{\label{Generalizacao4retas}Four straight lines}
\end{figure}

    (ii)  Given four straight lines $R_1$, $R_2$, $R_3$ and $R_4$, let $P_{12}$, $P_{13}$, $P_{14}$, $P_{23}$, $P_{24}$, $P_{34}$ be the intersection points accordingly to the indexes. In this way, our circles of interest will be those which contain points of the form $P_{ij}$, $P_{ik}$ and $P_{jk}$, with $i,j,k\in\{1,2,3,4\}$ and $i\neq j$, $i\neq k$ e $j\neq k$ (see Fig.\ref{Generalizacao4retas}).

    For example, the points $P_{12}$, $P_{13}$, $P_{23}$ lie on the circle $C_{123}$, such as the points $P_{ij}$, $P_{ik}$, $P_{jk}$ lie on the circle $C_{ijk}$. The order of the indexes is indifferent and in such a way we have formed four circles, $C_{123}$, $C_{124}$, $C_{134}$ and $C_{234}$. Then Theorem \ref{Miquel4} gives us the result expected, those four circles meet in the point $P_{1234}$.    \\

\begin{figure}[!h]
\center
  \includegraphics*[width=0.9\linewidth]{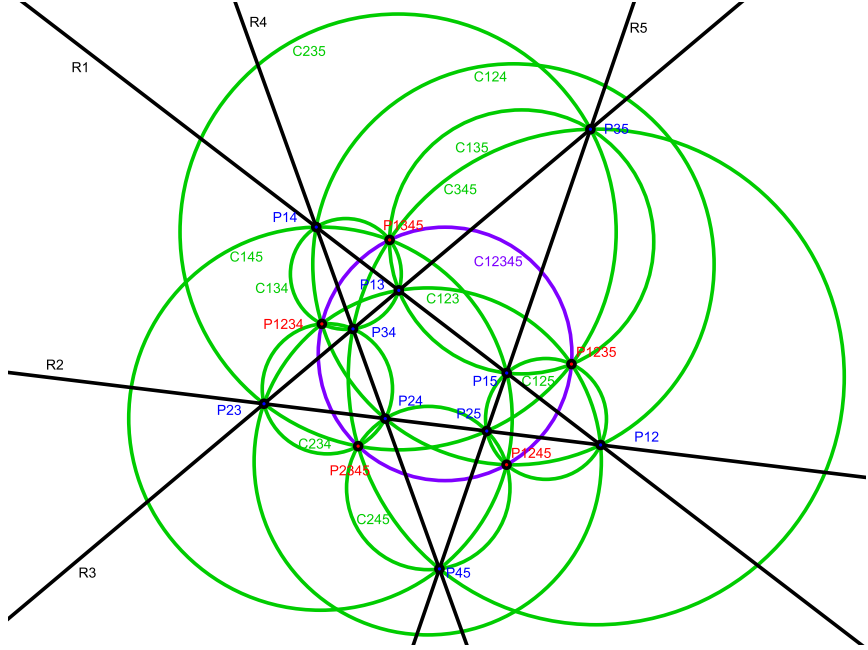}
  \caption{\label{Generalizacao5retas}Five straight lines}
\end{figure}

    (iii)  Given five straight lines $R_1$, $R_2$, $R_3$, $R_4$ and $R_5$, let $P_{12}$, $P_{13}$, $P_{14}$, $P_{15}$, $P_{23}$, $P_{24}$, $P_{25}$, $P_{34}$, $P_{35}$, $P_{45}$ be the intersection points accordingly to the indexes. Like in case (ii), the circles will be $C_{ijk}$, with $i,j,k\in\{1,\dots,5\}$, \emph{i.e.}, $C_{123}$, $C_{124}$, $C_{125}$, $C_{134}$, $C_{135}$, $C_{145}$, $C_{234}$, $C_{235}$, $C_{245}$ e $C_{345}$. Then case (ii) claims that a set of four of these circles meet in a point. For example, the circles $C_{123}$, $C_{124}$, $C_{134}$ and $C_{234}$ meet in the point $P_{1234}$. Analogously, we found other four points $P_{1235}$,  $P_{1245}$, $P_{1345}$ and  $P_{2345}$. Finally,  Theorem \ref{Miquel5} implies that these five points lie on the circle $C_{12345}$ (see Fig.\ref{Generalizacao5retas}).

    Let us prove the case (iii) again, now using the Lebesgue's idea. It is important to do so because this argument is indispensable for the last step of the induction.

    Consider the circles $C_{123}$, $C_{134}$, $C_{145}$ and $C_{125}$. Their intersections are: $C_{123}\cap C_{134}=\{P_{13},P_{1234}\}$, $C_{134}\cap C_{145}=\{P_{14},P_{1345}\}$, $C_{145}\cap C_{125}=\{P_{15},P_{1245}\}$ and $C_{125}\cap C_{123}=\{P_{12},P_{1235}\}$. Since $P_{13}$, $P_{14}$, $P_{15}$ and $P_{12}$, lie on the straight line $R_1$, accordingly by Lemma \ref{Lebesgue2}, the points $P_{1234}$, $P_{1345}$, $P_{1245}$ and $P_{1235}$ lie on the same circle.

    Note that we took the points which lay on the straight line $R_1$. Let us take now those which lie on $R_2$, that is, $C_{123}$, $C_{125}$, $C_{245}$ and $C_{234}$. Their intersections are: $C_{123}\cap C_{125}=\{P_{12},P_{1235}\}$, $C_{125}\cap C_{245}=\{P_{25},P_{1245}\}$, $C_{245}\cap C_{234}=\{P_{24},P_{2345}\}$ and $C_{234}\cap C_{123}=\{P_{23},P_{1234}\}$. Since $P_{12}$, $P_{25}$, $P_{24}$ and $P_{23}$, lie on the straight line $R_2$, accordingly by Lemma \ref{Lebesgue2}, the points $P_{1235}$, $P_{1245}$, $P_{2345}$ and $P_{1234}$ lie on the same circle.

    Nevertheless, observe that we have found two circles which have three points in common, so they are the same circle such as $C_{12345}$.\\

    (iv)  Given six straight lines $R_i$, with $i\in\{1,\dots,6\}$, we define their intersections $\{P_{ij}\}=R_i\cap R_j$ with $i\neq j$ and $i,j\in\{1,\dots,6\}$. Given these points we find the circles $C_{ijk}$ which contain $P_{ij}$, $P_{ik}$ and $P_{jk}$. Each set of four of these circles meet in a point $\{P_{ijkl}\}=C_{ijk}\cap C_{ijl}\cap C_{ikl}\cap C_{jkl}$, where $i$, $j$, $k$, $l$ are distinct values of the set $\{1,\dots,6\}$. Then we can determine six circles $C_{ijklm}$. Finally, we want to prove that these six circles meet in the point  $P_{123456}$.

    Let us secure indexes 5 and 6 and take the circles $C_{56123}$, $C_{563}$, $C_{564}$ and $C_{56124}$. So we have the following intersections: $C_{56123}\cap C_{563}=\{P_{1356},P_{2356}\}$, $C_{563}\cap C_{564}=\{P_{56},P_{3456}\}$, $C_{564}\cap C_{56124}=\{P_{1456},P_{2456}\}$ and $C_{56124}\cap C_{56123}=\{P_{1256},A\}$, where $A$ is the unknown point. By case (iii), we know that the points $P_{1356}$, $P_{56}$, $P_{1456}$ and $P_{1256}$ lie on the circle $C_{156}$, and then, by Lemma \ref{Lebesgue}, it follows that the points $P_{2356}$, $P_{3456}$, $P_{2456}$ and $A$ lie on one circle. Since $C_{23456}$ contains these three points, then $A\in C_{23456}$. On the other hand, we know that $P_{2356}$, $P_{56}$, $P_{2456}$ and $P_{1256}$ lie on the circle $C_{256}$, so $P_{1356}$, $P_{3456}$, $P_{1456}$ and $A$ lie on one circle. However, $C_{13456}$ contains these three points, that is, $A\in C_{13456}$.

      So we have proved that $A\in C_{12356}$, $A\in C_{12456}$, $A\in C_{23456}$ and $A\in C_{13456}$.

     Now let us secure indexes 1 and 2 and take the circles $C_{12356}$, $C_{123}$, $C_{124}$ and $C_{12456}$. So we have the following intersections: $C_{12356}\cap C_{123}=\{P_{1235},P_{1236}\}$, $C_{123}\cap C_{124}=\{P_{12},P_{1234}\}$, $C_{124}\cap C_{12456}=\{P_{1245},P_{1246}\}$ and $C_{12456}\cap C_{12356}=\{P_{1256},A\}$. Analogously, we know that $P_{1235}$, $P_{12}$, $P_{1245}$ and $P_{1256}$ lie on the circle $C_{125}$, so by Lemma \ref{Lebesgue}, it follows that the points $P_{1236}$, $P_{1234}$, $P_{1246}$ and $A$ lie on one circle. Since $C_{12346}$ contains these three points, then $A\in C_{12346}$. On the other hand, we know that $P_{1236}$, $P_{12}$, $P_{1246}$ and $P_{1256}$ lie on the circle $C_{126}$, so $P_{1235}$, $P_{1234}$, $P_{1245}$ and $A$ lie on one circle. However, $C_{12345}$ contains these three points, that is, $A\in C_{12345}$.

    We have proved that $A$ is the common point of the six circles, that is, $A=P_{123456}$.\\

\underline{\textbf{Induction hypothesis:}}

Suppose $n\geq2$:

\begin{itemize}
\item
    Given $2n$ straight lines, there are $2n$ circles meeting in the point $P_{12\dots (2n)}$;
\item
    Given $2n+1$ straight lines, there are $2n+1$ points lying on the circle $C_{12\dots(2n+1)}$.
\end{itemize}
${}$

\underline{\textbf{Induction step:}}

    We would like to prove that these statements are valid for $n+1$, \emph{i.e.}, we must regard two different cases, $2n+2$ and $2n+3$ straight lines.

    In order to make our notation simpler, let $P^*_i$ be the point which does not contain the index $i$, that is, $P^*_i=P_{12\dots (i-1)(i+1)\dots n}$, depending on the value of $n$ naturally. The same applies to circles.

    Let us take $2n+2$ straight lines and let us prove that the $2n+2$ determined circles intersect in the point $P_{12\dots (2n+2)}$.

    For each index $i=3$, 4, $\dots$, $2n$, $2n+1$, take the circles $C^*_2$, $C^*_{2i(i+1)}$, $C^*_{1i(i+1)}$ and $C^*_1$. Then,
    $$C^*_2\cap C^*_{2i(i+1)}=\{P^*_{2i},P^*_{2(i+1)}\}$$
    $$C^*_{2i(i+1)}\cap C^*_{1i(i+1)}=\{P^*_{12i(i+1)},P^*_{i(i+1)}\}$$
    $$C^*_{1i(i+1)}\cap C^*_1=\{P^*_{1i},P^*_{1(i+1)}\}$$
    $$C^*_1\cap C^*_2=\{P^*_{12},A\}$$

    Note that the intersection points have exactly $2n$ or $2n-2$ indexes, except for point $A$. By the induction hypothesis, they lie on a circle. The points $P^*_{2i}$, $P^*_{12i(i+1)}$, $P^*_{1i}$, $P^*_{12}$ lie on the circle $C^*_{12i}$. So Lemma \ref{Lebesgue} implies that the points $P^*_{2(i+1)}$, $P^*_{i(i+1)}$, $P^*_{1(i+1)}$ and $A$ lie on a circle. However, $C^*_{i+1}$ contains these three points and, consequently, we must have $A\in C^*_{i+1}$.

    On the other hand, the points $P^*_{2(i+1)}$, $P^*_{12i(i+1)}$, $P^*_{1(i+1)}$, $P^*_{12}$ lie on $C^*_{12(i+1)}$. So Lemma \ref{Lebesgue} implies that the points $P^*_{2i}$, $P^*_{i(i+1)}$, $P^*_{1i}$ and $A$ lie on a circle. Nevertheless, $C^*_i$ contains these three points and, consequently, we must have $A\in C^*_i$.

    Observe that $A\in C^*_1$ and $A\in C^*_2$, by construction. Therefore we have proved that $A\in C^*_i$ for $i=1$, 2, $\dots$, $2n+1$, $2n+2$. Thus, the $2n+2$ circles of interest meet on the point $A=P_{12\dots(2n+2)}$.

    Now let us take $2n+3$ straight lines and let us prove that the $2n+3$ determined points lie on the circle $C_{12\dots (2n+3)}$.

    For each $i=1$, 2, 3, $\dots$, $2n$, we fixate the indexes 1, 2, $\dots$, $i-1$, $i+4$, $\dots$ $2n+2$, $2n+3$ and take the circles $C^*_{(i+2)(i+3)}$, $C^*_{i(i+3)}$, $C^*_{i(i+1)}$, $C^*_{(i+1)(i+2)}$. Then,
    $$C^*_{(i+2)(i+3)}\cap C^*_{i(i+3)}=\{P^*_{i(i+2)(i+3)},P^*_{i+3}\}$$
    $$C^*_{i(i+3)}\cap C^*_{i(i+1)}=\{P^*_{i(i+1)(i+3)},P^*_i\}$$
    $$C^*_{i(i+1)}\cap C^*_{(i+1)(i+2)}=\{P^*_{i(i+1)(i+2)},P^*_{i+1}\}$$
    $$C^*_{(i+1)(i+2)}\cap C^*_{(i+2)(i+3)}=\{P^*_{(i+1)(i+2)(i+3)},P^*_{i+2}\}$$

    Note that the first point of each intersection set has exactly $2n$ indexes. So by induction hypothesis, the four intersection points lie on a circle. Lemma \ref{Lebesgue} implies that the points $P^*_i$, $P^*_{i+1}$, $P^*_{i+2}$ e $P^*_{i+3}$ lie on a circle too, for each $i=1,2,3,\dots,2n$.

    Besides that, all the circles so found are coincident, since each pair of these circles with $i=k$ and $i=k+1$ have three points in commom. It follows that the $2n+3$ points $P^*_i$ lie on the circle $C_{12\dots(2n+3)}$. \hfill $\blacksquare$\\

\section*{\textbf{Acknowledgement}}

    I am grateful to Professor Marcos Craizer for his helpful comments and his academic orientation. My thanks to CNPq for financial support during the preparation of this paper.

\bibliographystyle{amsplain}

\end{document}